\documentclass{amsproc}
\usepackage{eucal}
\newcommand{\B}{{\mathcal B}}
\newcommand{\C}{{\mathbb C}}
\newcommand{\cone}{{\bf C}}
\newcommand{\ch}{{\rm ch}}
\newcommand{\Exp}{{\rm Exp}}
\newcommand{\q}{{\bf q}}
\newcommand{\Hom}{{\rm Hom}}
\newcommand{\ind}{{\rm ind}}
\newcommand{\LL}{{\mathcal L}}
\newcommand{\MU}{{\bf MU}}
\newcommand{\w}{{\omega}}
\newcommand{\T}{{\mathbb T}}
\newcommand{\Q}{{\mathbb Q}}
\newcommand{\vol}{{\rm vol}}
\newcommand{\V}{{\bf V}}
\newcommand{\X}{{\bf X}}
\newcommand{\Z}{{\mathbb Z}}
\newsymbol \hslash 207D
\begin{document}

\title{Cobordism of symplectic manifolds and asymptotic expansions}
 
\author{Jack Morava}
\address{Department of Mathematics, Johns Hopkins University,
Baltimore, Maryland 21218}
\email{jack@math.jhu.edu}
\thanks{The author was supported in part by NSF grant 9504234.}
\date{1 June 1998}

\begin{abstract} The cobordism ring $\B_{*}$ of symplectic
manifolds defined by V.\ L.\ Ginzburg is isomorphic to the
Pontrjagin ring of complex-oriented manifolds with free circle
actions. This provides an interpretation of the formal group
law of complex cobordism as a composition-law on certain asymptotic
expansions. An appendix discusses some related questions 
about toric manifolds. \end{abstract}
 
\maketitle
 
\section*{\bf Introduction}
 
The theory of geometric quantization is concerned with smooth
manifolds $V$ bearing a complex line bundle $\LL$ with connection
$\nabla_{\LL}$, such that the Chern-Weil form $\w_{\LL}$ of the
connection is symplectic, i.\ e.\ a nondegenerate closed 
two-form representing an integral class in $H^{2}_{dR}(V)$. A
compatible almost-complex structure on a symplectic manifold is
an endomorphism $j$ with square $-1$ of the tangent bundle of
$V$, such that the bilinear form $g_{\LL}(-,-) = \w_{\LL}(j -,-)$
is a Riemannian metric; these structures exist because the unitary
group is a retract of the automorphism group of a real symplectic
vector space. For the purposes of this note, the quadruple
$\V := (V,\LL,\nabla_{\LL},j)$ will be called a prequantized
manifold, though this usage is not quite standard.

Analysis on this class of manifolds is very rich [cf.\ e.g.\
17-19] but I will be concerned here with the reflection in
topology of some aspects of that subject. I suggest, in
particular, that the formal group law on complex cobordism
[introduced by thirty years ago [20] by Novikov, which has led to
such spectacular advances in our understanding of homotopy theory
[22]], has a natural analytic interpretation in terms of
asymptotic expansions. The first section of this paper summarizes
some well-known facts about heat kernels, which provide a useful
class of examples of expansions; the second is an account of
the structure of the cobordism ring defined by V.L.\ Ginzburg
[10] for symplectic manifolds, and the third contains some 
closing remarks. I would like to thank V.L.\ Ginzburg 
for very interesting discussions about his work.

It is a pleasure to dedicate this paper to Professor Sergei
Petrovich Novikov, in view of his lifelong interest in the
frontier between geometry and analysis. This Web version differs
from the original, which appears in the Novikov anniversary 
volume of the Proceedings of the Steklov Institute, by the
addition of an appendix concerned with related questions about
toric manifolds. \bigskip

\section*{\bf \S 1 Topological invariants of prequantized
manifolds} 

A prequantized manifold, in the terminology adopted
here, has a preferred metric as part of its structure, and thus
permits the usual constructions of Riemannian geometry. The
twisted signature operator on $\Lambda(T^{*}V) \otimes \LL$
defined by a Dirac operator [in the sense of [3 \S 3.6], i.\ e.\
$d + \delta$ acting on forms graded by the Hodge operator]
associated to the bundle of Clifford algebras constructed from
the Riemannian metric is a particularly accessible example. The
Laplace-Beltrami operator $\Delta_{\LL}$ on sections of this
bundle [3 \S 2.4] has a well-behaved heat kernel on $V \times V$,
and when $V$ is closed and compact of dimension $2d$, the
restriction of that kernel to the diagonal possesses an
asymptotic expansion $$\exp(-t^{-1} \Delta_{\LL}) \sim t^{d}
\sum_{k \geq 0} t^{-i} k_{i}(\LL) \; d\vol_{\LL}$$ as $t \mapsto
\infty$; the coefficients $k_{i}(\LL)$ are homogeneous local
functionals of the metric and $\w_{\LL}$, of great analytic and
geometric interest [9 \S 4.1]. A theorem of McKean and Singer
implies that the index of $\Delta_{\LL}$ equals the integral over
$M$ of this differential form; after some simplification we
recover the Atiyah-Singer formula $$\ind \; \Delta_{\LL} =
(\ch(\LL) L(V))[V] \;,$$ where $L(V)$ is the Hirzebruch
$L$-genus, and $[V]$ is the orientation class of $V$.

Now if $\V = (V,\LL,\nabla_{\LL},j)$ is prequantized, and $n$ is
a positive integer, then $$ [n](\V) := (V,\LL^{\otimes
n},\nabla_{\LL^{\otimes n}},j)$$ is also a prequantized manifold.
The behavior of the analytic properties of prequantized manifolds
under this `Adams operation' is of considerable interest in the
theory of geometric quantization, where $n \mapsto \infty$ is
known as the `semiclassical limit', cf.\ [15 \S 34]. The Chern
class of $\LL^{\otimes n}$ is $n\w_{\LL}$, so $g_{\LL^{\otimes
n}} = ng_{\LL}$; the symbol $\sigma(\Delta_{\LL^{\otimes n}})$
thus equals $n^{-1} \sigma(\Delta_{\LL})$, and homogeneity
implies that $k_{i}(\LL^{\otimes n})$ is polynomial in $n$, while 
$d\vol_{\LL^{\otimes n}} = n^{d} d\vol_{\LL}$. It is a corollary 
that $\ind \; \Delta_{\LL^{\otimes n}}$ is polynomial of degree $d$ 
in $n$; it is an analogue in some sense of the Hilbert polynomial in
algebraic geometry. 

This invariant has some interesting properties: it is additive under
disjoint unions of manifolds, and it is multiplicative under
cartesian product; indeed this is true locally of the heat
kernels themselves [9 \S 1.8]. Since it is a characteristic
number, it is also a cobordism invariant, in a sense which will
be made precise below; and its relation to the theory of Feynman
path-integrals is rather well-understood [11 \S 3]. I have described this
example because it is typical of a large class of invariants 
constructible from heat kernel expansions; but recent striking work of 
K\'orpas and Uribe [17] suggests that the ${\rm Spin}^{c}$ Dirac 
operator has even more natural analytic properties.

The main technical result in this paper is a proof
that Ginzburg's symplectic cobordism is a (cocommutative) Hopf
algebra, isomorphic to the complex bordism ring $MU_{*}B\T$ of
manifolds with free circle actions; an immediate corollary is
that a certain class of additive functionals $\Phi$ on
prequantized manifolds (those such that $\Phi ([n](\V))$ possesses a
formal asymptotic expansion in $n$) forms a (cocommutative) Hopf
algebra, Cartier dual to Ginzburg's. This suggests as an
interesting question for the future, the possible existence of
some similar kind of algebraic structure on the heat
kernel invariants themselves, rather than on their global
integrals. 

\section*{\bf\S 2 Cobordism of symplectic and prequantized 
manifolds}
 
\noindent
{\bf 2.1} A (generalized) symplectic manifold is a pair $(V,\w)$,
with $V$ an oriented manifold and $\w$ a maximally nondegenerate
integral closed 2-form on $V$. This is the usual definition, if
$V$ is even-dimensional; if not, it entails that the kernel of
$\w$, viewed as a homomorphism from the tangent bundle of $M$ to
its cotangent bundle, has an everywhere one-dimensional kernel.
$\B_{*}$ denotes the graded ring of even-dimensional closed
compact symplectic manifolds, under the equivalence relation
defined by cobordism of such structures. 

Ginzburg [10 \S 1.7] defines a homomorphism $$\gamma : \B_{*}
\rightarrow  MU_{*}B\T$$ from his symplectic cobordism ring to
the complex bordism ring of the classifying space for complex
line bundles, which assigns to $(V,\w)$ a triple $(V,\LL,j)$,
consisting of a complex line bundle $\LL$ over $M$ with Chern
class represented by $\w$, together with an almost-complex
structure $j$ on the tangent bundle of $V$, such that the
bilinear form $$\w(j -,-)$$ is positive definite. He then shows
that this homomorphism is injective, and that it becomes an
isomorphism after tensoring with the rationals. It will be
convenient here to modify this construction very slightly, by
taking $B\T$ to be a classifying space for bundles with
connection [19]; the proof that $\gamma$ is an isomorphism then
implies that the forgetful map from the cobordism ring of
prequantized manifolds to $\B_{*}$ is an isomorphism as well.

The multiplicative monoid $\Z^{\times}$ of integers acts by ring
automorphisms of $\B_{*}$, sending $(V,\w)$ to $$[n](V,\w) =
(V,n\w) \;;$$ let $$\B^{0}_{*} = H^{0}(\Z^{\times},\B_{*})$$ be
the subring of invariant elements. There is a similar action of
$\Z^{\times}$ by ring automorphisms of $MU_{*}B\T$ induced from
the action $$[n] : z \mapsto z^{n}, z \in \T = \{ z \in \C| |z| =
1\}$$ on the circle group, or equivalently by the operation
$$[n](\LL) = \LL^{\otimes n}$$ on line bundles, and the
homomorphism $\gamma$ respects this action. The subring of
$\Z^{\times}$ invariants of $MU_{*}B\T$ is the cobordism ring
$MU_{*}$ of a point, so there is an induced monomorphism
$$\gamma^{0} : \B^{0}_{*} \rightarrow MU_{*} \;.$$ It is perhaps
surprising that this homomorphism has a relatively large image:
\bigskip

\noindent {\bf 2.2 Proposition}  {\it The cobordism classes $\C
P_{n} \in MU_{*}$ of the complex projective spaces lie in the
image of $\gamma^{0}$\;.} \medskip 

\noindent {\it Proof:}  We will construct polynomials in
symplectic manifolds of the form $(\C P_{m},k\w)$, where $n \geq
m$ and $k$ is a suitable positive integer, together with a
complex two-torus $X$ carrying a standard symplectic structure,
which map to the class of $\C P_{n} \in MU_{2n}$; here $\w$ is
the standard (Fubini-Study) form. The argument uses
characteristic numbers: if $(V,\w)$ is a compact closed
symplectic manifold, $$\sum_{I} (c_{I}\w^{n-|I|})[V]
t^{I}\w_{n-|I|}$$ will denote its image under the composition of
$\gamma$ with the Chern homomorphism $$MU_{*}B\T \rightarrow \Hom
(H^{*}\MU,H_{*}B\T) = \Z[t_{n},\w_{n}\;|\; n \geq 1]\;,$$ where
$\w_{n}$ is the $n$th divided power of $\w$, $I = i_{1},i_{2},
\dots$ is a multiindex of weight $|I| = \sum k i_{k}$, $c_{I}$ is
the polynomial in Chern classes of $V$ defined by the monomial
symmetric function $m_{I}$, and $t^{I} = \prod_{k \geq 1}
t_{k}^{i_{k}}$. For example, $$\begin{array}{lll} (\C P_{1},\w)
\longmapsto 2t_{1} + \w,\\ (\C P_{2},\w) \longmapsto 3t_{1}^{2} +
3t_{2} + 3t_{1}\w + \w_{2}, \\ (\C P_{3},\w) \longmapsto
4t_{1}^{3} + 12t_{1}t_{2} - 20t_{3} + (6t_{1}^{2} + 2t_{2})\w +
4t_{1}\w_{2} + \w_{3} \;, \end{array}$$ etc. The Chern
homomorphism is injective on the complex cobordism ring, and it
defines an identification of the polynomial $\Q$-algebra
generated by the $\C P_{n}$'s with the ring $\Q[t_{n}\;|\; n \geq
1]$. In particular, a polynomial in the $t$'s with rational
coefficients can be expressed uniquely as a sum of products $\C
P^{I} = \prod \C P_{k}^{i_{k}}$. We now show inductively that we
can construct a product of projective spaces, with standard
symplectic structures, such that all Chern numbers involving the
symplectic class vanish : if $\C P_{m}$ has been shown to lie in
$\B^{0}_{*}$ when $n>m$, we can write $$\mbox{(Chern polynomial
of)}\;\C P_{n} + \sum_{n>|I|} b_{I}\C P^{I} \w_{n-|I|} $$ for the
total symplectic Chern polynomial of $(\C P_{n},\w)$, with
coefficients $b_{I} \in \Q$. Let $k$ be the common denominator of
the collection $$\frac{b_{I}}{(n-|I|)!}$$ of rational numbers,
where $n > |I| ;$ then $$\C P_{n} = (\C P_{n},k\w) - \sum_{n>|I|}
k^{n-|I|-1} \frac{kb_{I}}{(n-|I|)!} \C P^{I}X^{n-|I|}$$ expresses
$\C P_{n}$ as a sum of elements of $\B_{*}$. For example,
$$\begin{array}{lll} \C P_{1} = (\C P_{1},\w) - X, \\ \C P_{2} =
(\C P_{2},2\w) - 3\C P_{1} X - 2X^{2},\\ \C P_{3} = (\C
P_{3},6\w) -(4\C P_{2} - 6\C P_{1}^{2})X - 36\C P_{1} X^{2} -
36X^{3}, \end{array}$$ etc. This completes the proof of the
proposition. \bigskip 

From now on, I will identify the class of the `classical'
projective space $\C P_{n} \in MU_{2n}$ with the corresponding
class in $\B^{0}_{2n}$; it should be distinguished from the class
of a `quantum' projective space $(\C P_{n},k\w)$. \bigskip

\noindent
{\bf 2.3} We now need some standard results about the complex
bordism and cobordism of $B\T$. If $(V,\LL)$ is a
complex-oriented manifold $V$ together with a complex line bundle
$\LL$ over it, then Quillen's Euler class [14] is the cobordism
element $$\q(\LL) = [s^{-1}(0_{\LL})] = V_{0} \in MU^{2}(V)$$
defined by the inverse image of the zero-section of $\LL$, under
a generic smooth section $s$. This defines a hyperplane section
homomorphism $$\q : MU_{*}B\T \rightarrow  MU_{*-2}$$ of
$MU_{*}$-modules, which can also be interpreted as the Kronecker
product of $(V,\LL)$ with the first Chern class in $MU^{*}B\T$;
more generally, there are homomorphisms $$\q^{i} : MU_{*}B\T
\rightarrow MU_{*-2i}$$ defined as Kronecker products with the
$i$th power of the Chern class; in geometric terms this is just
the $i$-fold transversal intersection $V_{0} \cap \dots \cap
V_{0}$ of hyperplane sections of $V$. Since $MU^{*}B\T$ and
$MU_{*}B\T$ are dual under the Kronecker pairing, there is a
basis $b_{k}$ of the latter module characterized by the formula
$$\q^{i}(b_{k}) = \delta_{i,k}\;.$$ The $i$th power of the Chern
class in $MU^{2i}B\T$ can be interpreted as the complex
projective subspace of $\C P_{\infty}$ of complex codimension
$i$. Let $$\beta_{k} : \C P_{k} \rightarrow \C P_{\infty}$$
denote the element of $MU_{2k}B\T$ defined by $(\C P_{k},\w)$;
then $$\q^{i}(\beta_{k}) = \C P_{k-i}\; ,$$ this being zero if
$i>k$. It follows that $$\beta_{i} =\sum_{k} \C P_{i-k} b_{k}
\;,$$ which implies that the classes $b_{i}$ lie in the subring
$\B_{*}$ of $MU_{*}B\T$: for $\beta_{i}$ certainly lies in
$\B_{*}$, while $\C P_{n} \in \B^{0}_{*}$ by the proposition
above. The assertion is thus a consequence of the invertibility
of the upper-triangular matrix $\C P_{i-j}$. 

The bordism ring of $B\T$ (with the Pontrjagin multiplication) is
known to be Cartier dual to $MU^{*}(B\T)$: if $$b(u) = \sum_{n
\geq 0} b_{n}u^{n} \in MU_{*}[[u]] \;,$$ then $$b(u)b(v) = b(u
+_{F} v) \leqno(*)$$ in $MU_{*}[[u,v]]$, where $u +_{F} v =
F(u,v)$ is the formal group sum in $MU^{*}[[u,v]]$, cf.\ [23 \S
3.4]. Now, however, we know that this equation holds in $\B_{*}$.
If we write $\beta(u)$ for the corresponding generating function
for the $\beta$'s, then $$b(u) = {\C P(u)}^{-1} \beta(u)\;,$$
with $\C P(u)$ the generating function for the classical
projective spaces. In this notation, the diagonal for the Hopf
algebra $\B_{*}$ is characterized by $$\C P(u) \Delta \beta(u) =
\beta(u) \otimes \beta(u)\;.$$

\noindent {\bf 2.4 Proposition} {\it The completed localization
$$\B_{*}[[\C P^{-1}_{1}]] \rightarrow (MU_{*}B\T)[[\C P^{-
1}_{1}]]$$ of $\gamma$ is an isomorphism.} \medskip

\noindent {\it Proof:} The element $$b_{1} = \beta_{1} - \C
P_{1} = X$$ is a unit in $\B_{*}[[\C P_{1}^{-1}]]$, so the power
series $b(u) - 1$ possesses a formal inverse $\tilde b(u)$ with
coefficients in $\B_{*}[[\C P^{-1}_{1}]]$. Formula $(*)$ can thus
be restated as $$\tilde b(b(u)b(v) - 1) = u +_{F} v \;;$$ it
follows that the coefficients of the formal group law $F(u,v)$,
which are a priori elements of $MU_{*}$, all lie in
$\B_{*}[[\C P^{-1}_{1}]]$. But the coefficients of the formal
group law [4] generate $MU_{*}$, so
$MU_{*}$ is a subring of $\B_{*}[[\C P^{-1}_{1}]]$. This latter
ring is then a subring of $(MU_{*}B\T)[[\C P^{-1}_{1}]]$ which
contains $MU_{*}[[\C P^{-1}_{1}]]$ as well as the classes $b_{n}$
(by the preceding paragraph); but the $b$'s are a basis for
$MU_{*}B\T$ over $MU_{*}$, so $\B_{*}[[\C P^{-1}_{1}]]$ contains
the full ring $(MU_{*}B\T)[[\C P^{-1}_{1}]]$, and is thus
isomorphic to it. \bigskip 

\noindent {\bf 2.5 Theorem} {\it The homomorphism $$\gamma:
\B_{*} \rightarrow MU_{*}B\T$$ is an isomorphism.} \medskip 

\noindent {\it Proof:}  Since $\C P_{1}$ is
$\Z^{\times}$-invariant, the completion $\B_{*}[[\C P^{-1}_{1}]]$
inherits a $\Z^{\times}$-action, with invariant subring
$\B^{0}_{*} [[\C P^{-1}_{1}]$] contained in $MU_{*}[[\C
P^{-1}_{1}]]$. The argument above shows that this inclusion is in
fact an isomorphism, and to complete the proof of the theorem it
is enough to show that $$\gamma^{0} : \B^{0}_{*} \rightarrow
MU_{*}$$ is an isomorphism.

Now consider the map induced by this homomorphism on
modules of indecomposables: let $I_{\B}$ be the (graded) ideal of
positive-dimensional elements of the domain, and $I_{MU}$ the
corresponding ideal in the range; then by Ginzburg's theorem the
induced homomorphism $$I_{\B}/I_{\B}^{2} \rightarrow  
I_{MU}/I_{MU}^{2}$$ becomes an isomorphism after tensoring with
the rationals. Since both of these quotient modules are free of
rank one in even dimensions, and zero otherwise, this quotient
homomorphism is specified by one nonzero integer (its value on a
generator) for each even dimension. On the other hand this map is
also an isomorphism after completing with respect to
$\C P^{-1}_{1}$; thus these integers are all units.\bigskip 

In the time since the first draft of this paper was written, an
independent proof for this theorem has been found by A.\ Baker
[2], based on symplectic structures on Milnor's hypersurface
generators. \bigskip

\section*{\bf \S 3 A Hopf algebra of functionals} 

The existence of a Hopf algebra structure on the
cobordism ring of symplectic manifolds is in some ways rather
surprising. This section is concerned with some consequences of
this result, expressed in terms not of symplectic manifolds 
but rather of functionals on them. Perhaps it is best
to start with a central example: the formal series $$\hslash(\V) =
\sum_{n \geq 1}\frac{\C P_{n-1}}{n} \cap^{n} V_{0}$$ [with $\cap^{n} V_{0}$
interpreted as the cobordism class of an $n$-fold transversal intersection]
defines a homomorphism $$\B_{*} \rightarrow \B^{0}_{*}
\otimes \Q $$ of degree two of graded abelian groups, such that \medskip

\noindent i) if $M \in \B^{0}_{*}$ is classical, and $\V \in \B_{*}$ 
is quantum, then $$\hslash(M\V) = M \hslash(\V) \; ,$$ while \medskip

\noindent ii) $\hslash(V,n\w) = n \hslash(V,\w)$ . \medskip

\noindent The first of these properties holds for many
interesting linear functionals on $\B_{*}$, for example the 
twisted signature considered in the introduction; but the second 
is rather special.
In fact it is a restatement of Quillen's theorem that the module
$\B^{*}$ of graded $\B^{0}_{*}$-linear homomorphisms from
$\B_{*}$ to $\B^{0}_{*}$ is a Hopf algebra isomorphic to
$\B^{0}_{*}[[\q]]$, with $$\q = \Exp (\hslash) \;,$$ where $\Exp$
is the formal series inverse to Mi{\v s}{\v c}enko's logarithm.
It is tempting to regard $n^{-1} \C P_{n-1}$ as a cyclic
quotient of $\C P_{n-1}$ defined by the shift
$[z_{0}:\dots:z_{n-1}] \mapsto [z_{n-1}:z_{0}: \dots ]$,
and to think of the formula for $\hslash$ as a correction to the hyperplane
section construction which keeps track of some kind of bubbling.

In general, any linear functional on cobordism classes of
prequantized manifolds which  satisfies condition i) factors
through $\B^{0}_{*}[[\q]]$; for example, the universal
$K$-theoretic characteristic number homomorphism [6] can be
characterized as the $\B^{0}_{*}$-linear map which sends $b_{1} =
X$ to a unit. \bigskip 

\noindent
{\bf Appendix: some questions about cobordism of toric manifolds}
\bigskip

\noindent
{\bf 1} Toric varieties are by now familiar objects in algebraic
geometry, but this appendix is concerned with variations on that
theme, and I will try to be careful about terminology. A toric
variety is a kind of orbifold, and hence has mild singularities,
but I will use the term toric manifold in the sense of Davis and
Januszkiewicz [7]; a smooth toric variety thus has an underlying
toric manifold, but toric manifolds form a slightly more general
class. \bigskip

\noindent
{\bf 2} By symplectic toric manifold I mean a compact $2n$-dimensional 
symplectic manifold with a Hamiltonian $\T^{n}$-action, which is 
moreover {\bf prequantizeable} in the sense of Guillemin [14 Th. 3.2]: 
such a thing is defined by a Delzant polytope [14 Th. 1.8] with 
vertices at integer lattice points, and the Chern class of its 
canonical line bundle is calculated in [14 appendix 2]. If we ignore 
the torus action, then the underlying symplectic manifold defines an 
element of Ginzburg's symplectic cobordism ring $\B_{*}$. [There is a 
small but potentially confusing point here, in that elements of 
Ginzburg's ring are equivalence classes of compact {\bf oriented} manifolds
endowed with a nondegenerate closed integral two-form $\omega$. I
will assume that the orientation is defined by $\omega^{n}$.]
\bigskip

\noindent
{\bf 3} Ginzburg's cobordism ring is canonically isomorphic to the ring
$MU_{*}B\T$ defined by complex-oriented manifolds together with a
complex line bundle; forgetting this line bundle defines a
homomorphism to the complex bordism ring $MU_{*}({\rm pt})$. A Delzant 
polytope thus defines a cobordism class, and Buchstaber
and Ray [5] have recently published an elegant proof that this
construction is surjective: {\bf every} complex-oriented manifold
is cobordant to a (sum of) toric manifold(s); at least, if we 
allow polytopes of dimension zero. In fact Conner and Floyd
constructed multiplicative generators for the cobordism ring (at
least, for each prime $p$) as Bott towers with obvious toric
manifold structures [6 \S 42]; but this happened a few years before the
notion of toric variety appeared in print. It seems pretty clear
that the Cartesian product of toric manifolds is again a toric 
manifold; but it doesn't seem to me so immediate that the Cartesian 
product of two Delzant polytopes is another such. I don't know a 
reference for this presumed fact \dots  \bigskip

\noindent
{\bf 4} In light of these observations, it is natural to ask which elements 
of $\B_{*}$ `come from' symplectic toric actions; but we need to
be precise about how this is to be interpreted. If $(V_{\Delta},\omega)$
is the class defined by an integral Delzant polytope $\Delta$,
then it seems reasonable to think of the elements of the ray 
$$\{ (V_{\Delta},n\omega) \in MU_{*}(B\T) \;|\; n \in \Z_{>0}\}$$
as arising in this way. The disjoint union of two
symplectic manifolds defined by torus actions ought also to be in
this class. There is thus a kind of {\bf effective cone} $\cone$ in
$MU_{*}(B\T)$, generated by sums of elements of such rays; this
cone will be closed under multiplication, if the remark above is
not a mistake. It would be very interesting to know more about $\bf C$. 
\bigskip

\noindent
I want to thank Ginzburg, Guillemin, and Karshon for very helpful
correspondence about this question. In response to some very naive
conjectures, they observed that in dimension two, $\B_{*}$ has rank two, 
with the projective line $\C P_{1}$ and the symplectic torus $\X$ 
as generators; the latter cannot be toric. In fact it seems plausible
that no class in principal ideal $(\X)$ can be toric; if so, then 
the quotient of the localization ${\bf C}^{-1}\B_*$ by the image of $(\X)$ 
will be the quotient field of $MU_{*}$, and the localization itself 
would be a local ring. In some ways this ring is a natural context for
questions about geometric localization theorems along the lines of [13].
\bigskip

\noindent
{\bf 5} A related circle of questions concerns the relation of toric 
manifolds to the `normally split' cobordism theory $ML$ defined by
manifolds endowed with a stable splitting of the normal bundle as a
sum of line bundles [1]. Let $P$ be a simple convex polytope of dimension 
$n$, with $F$ its set of codimension one faces, and let $\lambda$ be a 
characteristic function [i.e.\ a suitable linear function from the free 
$\Z$-module generated by $F$, to the free $\Z$-module on standard generators 
$e_{i}$]. The Stanley-Reisner {\bf face ring} of $P$ is a 
(combinatorially-defined) quotient of the polynomial ring on generators 
$x_{f}$ indexed by the faces of $P$; let $$\lambda^{*} : H^{*}_{\T^{n}}
({\rm pt}) := \Z[e_{1},\dots,e_{n}] \rightarrow \Z[x_{f}|f \in F]$$ be 
the ring homomorphism induced by the dual of $\lambda$. Then there is a 
manifold $M := M(P)_{\lambda}$ with $\T^{n}$ action, such that the 
equivariant cohomology $H^{*}_{\T^{n}}(M)$ is isomorphic, as (free [7 \S 
4.12]) $H^{*}_{\T^{n}}({\rm pt})$-algebra, to the face ring of $P$ 
[7 \S 4.8]; moreover, the (equivariant) tangent bundle of $M$ splits stably, 
as a sum of complex line bundles indexed by $F$ [7 \S 6.6], such that the 
$f$th line bundle has Chern class $x_{f}$. \bigskip

\noindent
We can thus define an {\bf exhaustive} characteristic class for such a toric 
manifold as the product $$\prod_{f \in F} (1 - x_{f}z_{f})^{-1} \in 
H^{*}_{\T^{n}}(M)[z_{f}|f \in F] \; ,$$ and we can similarly define an 
exhaustive characteristic number polynomial of $M$ as the homogeneous 
component of degree $2n$ of the image of this polynomial, in the quotient of 
the equivariant cohomology by the ideal $(e_{1}, \dots , e_{n})$. Since the
equivariant cohomology is free over $H^{*}_{T^{n}}({\rm pt})$, this quotient 
can be identified with the ordinary (nonequivariant) cohomology of $M$, 
which is isomorphic to $\Z$ in degree $2n$; the polynomial is thus a 
homogeneous element of $\Z[z_{f}|f \in F]$, of (homological) degree $2n$.
This suggests the question: does a toric manifold define, in some natural 
way, an $ML$-cobordism class?

\bibliographystyle{amsplain}

\end{document}